\documentclass[12pt]{article}
\usepackage{amsfonts,amssymb,amscd,amsthm}

\def\be{\begin{equation}}
\def\ee{\end{equation}}
\def\bea{\begin{eqnarray}}
\def\eea{\end{eqnarray}}

\def\z{{z}}
\def\h{z}
\def\anticonm#1#2{\left [ {#1},{#2} \right ]_+ }

\def\pois#1#2{\left\{ {#1},{#2} \right\}}         
         
\def\conm#1#2{\left [ {#1},{#2} \right ]}

\begin{document}


\thispagestyle{empty}
\hfill \today

\vspace{2.5cm}

\begin{center}
\bf{\LARGE Poisson-Hopf limit of 
quantum algebras}
\end{center}

\bigskip\bigskip

\begin{center}
A. Ballesteros$^1$, E. Celeghini$^2$  and M.A. del Olmo$^3$
\end{center}

\begin{center}
$^1${\sl Departamento de F\'{\i}sica, Universidad de Burgos, \\
E-09006, Burgos, Spain.}\\
\medskip

$^2${\sl Departimento di Fisica, Universit\`a  di Firenze and
INFN--Sezione di
Firenze \\
I50019 Sesto Fiorentino,  Firenze, Italy}\\
\medskip

$^3${\sl Departamento de F\'{\i}sica Te\'orica, Universidad de
Valladolid, \\
E-47005, Valladolid, Spain.}\\
\medskip

{e-mail: angelb@ubu.es, celeghini@fi.infn.it, olmo@fta.uva.es}

\end{center}

\bigskip

\bigskip

\begin{abstract}
The Poisson-Hopf analogue of  an arbitrary quantum algebra $U_z(g)$ is constructed by introducing
a one-parameter family of quantizations $U_{z,\hbar}(g)$ depending explicitly on $\hbar$ and by taking the appropriate $\hbar\to 0$ limit.  The $q$-Poisson analogues of the $su(2)$  algebra are discussed and the novel $su_q^{\cal P}(3)$ case is introduced. The $q$--Serre relations are also extended to the Poisson limit. This  approach opens the perspective for possible applications of higher rank $q$-deformed Hopf algebras in semiclassical contexts.
\end{abstract}

\vskip 1cm

MSC: 17B63, 17B37, 81R50   

\vskip 0.4cm

\noindent Keywords: Poisson algebras, semiclassical limit, deformations, Hopf algebras, quantum groups,  $su(3)$ algebra

\vfill\eject


\section{Introduction\label{introduccion}}

Quantum groups were initially introduced as quantizations of Poisson-Lie groups associated 
to certain solutions of the classical Yang-Baxter equation. In this context, the deformation 
parameter were taken as 
$q=e^{\h}$, 
where $\h$ is the constant that governs the 
noncommutativity of the algebra of observables given by the quantum group entries, and quantum algebras were 
obtained as the Hopf algebra dual of quantum groups (for a detailed discussion, see \cite{KR}-\cite{majid} and references therein). In the case of the 
transition from classical to quantum physical models, the deformation parameter $\h$ was interpreted as the 
Planck constant $\hbar$.

However, in more general contexts $\h$ is a parameter whose geometric/physical meaning has to be elucidated for each particular case. 
In fact, quantum groups and quantum algebras were soon considered as `abstract'  Hopf algebras 
(being both noncommutative and noncocommutative) in order to explore whether these new 
objects can be considered as new symmetries of some physically relevant systems. The 
keystone of this approach was the discovery of the $su_q(2)$ invariance of the Heisenberg 
spin XXZ chain \cite{PS,KulSkly}, that was followed by a number of results exploiting quantum 
algebra symmetries in two-dimensional models \cite{Sierra}. Indeed,  in  the XXZ chain the `quantum' deformation 
parameter $q$ is neatly identified with the anisotropy of the chain, which is completely 
independent of the (truly quantum) $\hbar$ constant. And the same independence with respect 
to $\hbar$ can be traced in many other physical applications of quantum algebras and groups, 
like for instance lattice systems  (where $q$ is related to the lattice 
length) \cite{phonons,magnons}, deformations of kinematical symmetries   (in which the 
deformation parameter is a fundamental scale, see \cite{Luki}--\cite{Afin} and references therein), effective nuclear models 
\cite{nuclei,Bonatsos}, etc.

The aim of this paper is to introduce a mathematical framework for quantum algebras in which 
both the deformation parameter $z = \log q$ and the Planck constant 
$\hbar$ are independently 
and simultaneously considered.  
This approach is described in section 2, where for any quantum algebra $U_z(g)$  
we construct  a one-parameter family of equivalent quantizations 
$U_{z,\hbar}(g)$ that depends explicitly on both $\hbar$ and $z$. Then, a 
$q$-Poisson Hopf algebra is 
obtained as the $\hbar\to 0$ limit, thus obtaining -as in the Lie case- a proper  classical-mechanical limit of 
quantum systems endowed with an arbitrary quantum algebra symmetry. This approach can be interesting to construct new (classical) integrable systems having a parameter ($z$) that can control their dynamical behaviour.

In section 3 the method is illustrated by applying it to the quantum deformations of $su(2)$, considering both the standard and the non-standard ones. In this elementary case the $q$-Poisson Hopf algebras obtained as the $\hbar\to 0$ limit look to be formally identical to the original quantum algebras. 
The case of $su(3)$, fully described in section 4,  is the first non-trivial one since by starting from the quantum algebra $U_{z,\hbar}(su(3))$ we obtain a $q$-Poisson $su(3)$ algebra (that we shall call $su_q^{\cal P}(3)$) which is quite different from the former. 
To the best of our knowledge, this is a new $q$-Poisson Hopf algebra that could be used, for instance, in order to construct integrable deformations of  higher rank classical Gaudin models \cite{Gau}-\cite{Lerma} by using the approach presented in \cite{BR} or to consider the semiclassical limit of the $U_{z}(su(3))$ dynamics from the viewpoint of \cite{DaP}. With this in mind, the explicit form of the two Casimir functions for $su_q^{\cal P}(3)$ are explicitly found. A concluding section closes the paper, in which the Poisson analogue of the $q$-Serre relations for higher rank $q$-Poisson Hopf algebras is consistently defined.


\section{Quantum algebras and Poisson-Hopf limit\label{Poissonlimit}}

Let us recall that a Lie bialgebra $(g,\delta)$ is a Lie algebra $g$
\be\label{constantesestructura}
[X_i,X_j]= f^k_{ij}\,X_k
\ee
together with a cocommutator map $\delta: g\rightarrow g\otimes g$ given by
\be
\delta(X_i)= c_i^{jk}\,X_j\otimes X_k
\label{tensorc}
\ee
such that $c_i^{jk}$ defines a (dual) Lie algebra and fulfills 
the appropriate compatibility condition (see \cite{CP,majid} for details). 

Given an arbitrary  Lie bialgebra $(g,\delta)$ the quantum algebra $(U_{z}(g),\Delta_z)$ (that is obtained through the analytic procedure described in \cite{ballesterosJCS08,ballesterosJPA08}) is 
the Hopf algebra deformation $(U_{z}(g),\Delta_z)$  of the universal enveloping algebra of $g$, 
$U(g)$, compatible with the deformed coproduct $\Delta_z(X)$ whose leading order terms are
$$
\Delta_z(X)= \Delta_0(X) + z\, \delta(X) + o[z^2].
$$

Let us now consider the  one-parameter family of equivalent Lie bialgebras $(g_\hbar,\delta)$ defined by (\ref{tensorc}) and
\be
[X_i,X_j]=\hbar\, f^k_{ij}\,X_k .
\label{hconm}
\ee
Note that when $\hbar=1$ we recover the original commutation relations (\ref{constantesestructura}).
If we quantize $(g_\hbar,\delta)$ by using the method described in~\cite{ballesterosJPA08} we obtain the two-parameter quantum algebra  $(U_{z,\hbar}(g),\Delta_{z,\hbar})$,  that depends explicitly on $\hbar$. Then, the Poisson limit ($\hbar \to 0$) of  $(U_{z,\hbar}(g),\Delta_{z,\hbar})$ can be uniquely defined and gives the Poisson-Hopf algebra $(Fun(g_z),\Delta_z^{\cal P})$. This $q$-Poisson algebra is just a Poisson-Lie structure on the group $g_z$ whose Lie algebra is determined by the dual $\delta^\ast$ of the Lie bialgebra map (\ref{tensorc}), {\em i.e.}, by the structure tensor $c_i^{jk}$ (see \cite{CP,majid}). The Poisson bracket on $Fun(g_z)$ is given by
\be\label{poislimit}
\pois{X}{Y}:=\lim_{\hbar\to 0}{\frac{\conm{X}{Y}}{\hbar}}
\ee
and the coproduct map
\be\label{colimit}
\Delta_z^{\cal P}(X):=\lim_{\hbar\to 0}{\Delta_{z,\hbar}}(X)
\ee
is a Poisson algebra homomorphism between $Fun(g_z)$ and $Fun(g_z)\otimes Fun(g_z)$.

If we deal with non-deformed Lie algebras, its coproduct is the primitive one and 
the commutation rules are linear: thus the Poisson limit (\ref{poislimit}) and (\ref{colimit}) leads to the same formal structure where commutation rules have been just replaced by Poisson brackets. On the contrary, quantum algebras introduce nonlinear functions of the generators both at the level of the commutation rules and of the coproduct. This implies that the $q$-Poisson structure given by the limit $\hbar\to 0$ can be formally different to the original quantum algebra structure. In fact, such limit  allows us to remove contributions in the deformation that arise as reordering terms, as shown in the following sections. 

Summarizing, in this paper we introduce  the following commutative diagram
\[ \begin{CD}
(U(g),\Delta_{0}) @<{\rm \hskip0.42cm \hbar=1 \hskip0.42cm}<<
(U_{\hbar}(g),\Delta_{0}) @>{\rm \hskip0.42cm \hbar\to 0 \hskip0.42cm}>> (\mbox{Fun}(g),\Delta_0^{\cal P})  \\
@AA{z\to 0}A @ AA{z\to 0}A @ AA{z\to 0}A  \\ 
(U_{z}(g),\Delta_{z}) @<{\rm \hskip0.42cm \hbar=1 \hskip0.42cm}<<(U_{z,\hbar}(g),\Delta_{z,\hbar}) @>{\rm \hskip0.42cm \hbar\to 0 \hskip0.42cm}>> (\mbox{Fun}(g_z),\Delta_z^{\cal P})  
\end{CD}\]
where  we focus on the lower right  corner,  $(\mbox{Fun}(g_z),\Delta_z^{\cal P})$, that we define in a constructive way by starting from any Lie bialgebra. 
 This general approach is illustrated in the following sections through the 
$su_q(2)$ and $su_q(3)$ examples. In particular, the $q$-Poisson algebra presented in section 4 is, to the best of our knowledge, the first example of a Hopf algebra deformation of the Poisson $su(3)$ algebra.  


\section{$q$-Poisson-Hopf algebras related to $su(2)$ }

In the $su(2)$ case,
two  well known quantum deformations do exist:  the standard one 
 \cite{Dri87,Jimbo} and  the non-standard
(or Jordanian) deformation (see \cite{Ohn}). As we shall see in the sequel, their corresponding $q$-Poisson Hopf algebras are formally equivalent to the quantum algebras from which they have been obtained.

\subsection{Standard $q$-Poisson algebra  $su_{q}^{\cal P}(2)$}

The $su(2)$ commutation rules are 
\be
[F_{12},F_{21}]= 2\,H \qquad [H,F_{12}]=F_{12} \qquad [H,F_{21}]=-F_{21}
\label{su2}
\ee
where $H=(H_1-H_2)/2$.
The standard $su(2)$ Lie bialgebra is given by 
\be
\delta (H)=0 \qquad
\delta (F_{12})=  H \wedge F_ {12} \qquad
\delta (F_{21})=  H \wedge F_ {21} .
\label{su2d}
\ee
The well-known quantum algebra deformation of (\ref{su2}) and (\ref{su2d}) reads
\be\label{suz}
[F_{12},F_{21}]= \frac{{\sinh(2\,z\,H)}}{z}\qquad
[H,F_{12}]=\,F_{12}
\qquad [H,F_{21}]=-\,F_{21}
\ee
\bea\label{cosu2}
&& \Delta_z (H) = H \otimes 1 + 1 \otimes H \nonumber\\[0.3cm]
&& \Delta_z (F_{12}) = e^{z\,H} \otimes F_{12} + F_{12} \otimes
e^{-z\,H} \\[0.3cm]
&& \Delta_z (F_{21}) = e^{z\,H} \otimes F_{21} + F_{21} \otimes
e^{-z\,H} 
\nonumber
\eea

Now if we consider the $\hbar$-parameter Lie algebra 
\be
[F_{12},F_{21}]= 2\,\hbar\,H \qquad [H,F_{12}]=\hbar\,F_{12} \qquad
[H,F_{21}]=-\hbar\,F_{21}
\label{su2h}
\ee
we would obtain a quantization in which the coproduct (\ref{cosu2}) does not formally change, but we have
\be
[F_{12},F_{21}]=\hbar\, \frac{\sinh(2\,z\,H)}{z}\,\qquad
[H,F_{12}]=\hbar\,F_{12}
\qquad [H,F_{21}]=-\hbar\,F_{21}
\ee
As a consequence, the $\hbar$-deformed Casimir operator is shown to be 
\be\begin{array}{lll}\label{casimirzh}
\ \hskip-0.5cm C_q &=&\displaystyle  \frac{{\sinh^2(z H)}}{z^2}
\cosh (z\hbar) +
\frac12 \;[F_{12}\, ,F_{21}]_+  
\\[0.35cm]
&=&\displaystyle \frac{{\sinh(z H)}}{z}   \; 
\frac{\sinh z( H + \hbar)}{z} + F_{21}\,F_{12}
= \displaystyle \frac{{\sinh(z H)}}{z}
 \frac{\sinh z(H - \hbar)}{z}+ F_{12}\,F_{21}
\end{array}
\ee
where the role of $\hbar$ can be easily appreciated.

We stress that  if we perform the following substitution
\[
H\rightarrow \hbar\,H
\qquad
z\rightarrow\,z/\hbar
\qquad
F_{12}\rightarrow \hbar\,F_{12}
\qquad
F_{21}\rightarrow \hbar\,F_{21}
\]
the parameter $\hbar$ can be reabsorbed 
and we get
 just the standard expressions  
(\ref{suz}) and (\ref{cosu2}), thus
confirming that the quantum algebra $U_z(su(2))$ depends on only one essential parameter.
However, since the limit $\hbar\to 0$ does not commute with this mapping, 
 the $\hbar$ parameter is essential in order to define the proper
$q$-Poisson algebra.

As a result,  the $q$-Poisson algebra  $su_{q}^{\cal P}(2)\equiv (\mbox{Fun}(g_z),\Delta_z^{\cal P})$ is obtained as the commutative algebra of functions endowed with the Poisson bracket and coproduct map coming, respectively, from the limits (\ref{poislimit}) and (\ref{colimit}). Namely,
\be
\pois{F_{12}}{F_{21}}=  \frac{\sinh(2\,z\,H)}{z} \qquad
\pois{H}{F_{12}}=F_{12}
\qquad \pois{H}{F_{21}}=- F_{21}
\ee
\be\begin{array}{l}
 \Delta_z^{\cal P} (H) = H \otimes 1 + 1 \otimes H\\[0.3cm]
 \Delta_z^{\cal P}  (F_{12}) = e^{z\,H} \otimes F_{12} + F_{12} \otimes
e^{-z\,H} \\[0.3cm]
 \Delta_z^{\cal P}  (F_{21}) = e^{z\,H} \otimes F_{21} + F_{21} \otimes
e^{-z\,H} 
\end{array}\ee
which is indeed a Poisson-Hopf algebra, as it can be checked by direct
computation. 
Finally, a unique $q$-deformed Casimir function is obtained from (\ref{casimirzh}) as
\[
C^{\cal P} _q = \lim_{\hbar\to 0}{C_q}=
\frac{1}{z^2} \; {\sinh^2(z\,H)} + F_{12}\,F_{21}.
\] 

\subsection{Non-standard $q$-Poisson   $su_{ns}^{\cal P}(2)$}

The same procedure can be applied to the non-standard deformation of $su(2)$. Explicitly, 
the non-standard  $\hbar$-parameter Lie bialgebra would be given by (\ref{su2h}) and
\be
\delta(F_{12})=0,\qquad \delta(H)=H\wedge F_{12},\qquad 
\delta(F_{21})=F_{21}\wedge F_{12} .
\label{ac}
\ee
The corresponding quantum algebra is \cite{Ohn} 
\bea
&& [H,F_{12}]=\hbar\,\frac{\sinh ( z F_{12}  )} z \qquad\qquad
[F_{12}  ,F_{21} ]= 2\,\hbar\,H\nonumber\\
&& [H,F_{21}]=-\frac{\hbar}{2}\left( F_{21}\,\cosh (z F_{12}  ) +\cosh (z F_{12} 
)
\,F_{21} \right)
\nonumber
\eea
\bea
&&  \Delta_{z,\hbar} (F_{12})  =1 \otimes F_{12}  + F_{12} \otimes 1\nonumber\\[0.3cm]
&& \Delta_{z,\hbar} (H) =e^{-z F_{12} } \otimes H + H\otimes
e^{z F_{12} } \nonumber \\[0.3cm]
&& \Delta_{z,\hbar} (F_{21}) =e^{-z F_{12} } \otimes F_{21} + F_{21}\otimes
e^{z F_{12} }
\nonumber
\eea
whose Casimir operator reads
\be
{\cal C}_{ns}=H^2  + \frac12\left(  \frac{\sinh z\,F_{12}}{z} \,F_{21} +
F_{21}\,\frac{\sinh z\,F_{12}}{z}\right)
+\frac{\hbar^2}{4}\cosh^2(z\,F_{12}).
\label{gc}
\ee

Now, by performing the same $\hbar\to 0$ limit given by (\ref{poislimit}) and (\ref{colimit}) we get
the non-standard $q$-Poisson-Hopf algebra $su_{ns}^{\cal P}(2)$
\bea
&& \pois{H}{F_{12}}=\frac {\sinh ( z F_{12}  )}{z} \quad\qquad 
\pois{F_{12}}{F_{21} }= 2\,H\nonumber\\
&& \pois{H}{F_{21}}=-F_{21}\,\cosh (z F_{12})  \quad
\label{ae} 
\eea
\bea
&&  \Delta_z^{\cal P} (F_{12})  =1 \otimes F_{12}  + F_{12} \otimes 1\nonumber\\[0.3cm]
&& \Delta_z^{\cal P} (H) =e^{-z F_{12} } \otimes H + H\otimes
e^{z F_{12} } \label{ad}\\[0.3cm]
&& \Delta_z^{\cal P} (F_{21}) =e^{-z F_{12} } \otimes F_{21} + F_{21}\otimes
e^{z F_{12} }
\nonumber
\eea
And the
Casimir function is
\be
{\cal C}^{\cal P}_{ns}=H^2  + \frac {\sinh (z\,F_{12})}{z}\,F_{21}
\label{gc}
\ee


\section{The $q$-Poisson $su_q^{\cal P}(3)$ algebra} 

In this case we will refer to \cite{ballesteros06,ballesteros07}, where the standard deformation of  $su(3)$ is obtained by starting
from the Weyl-Drinfeld basis of the bialgebra where all roots are
well defined. In this way a
complete description of the whole structure for $u_q(3) \equiv
su_q(3)\oplus u(1)$, real form of $A_2^q\oplus A_1$, is obtained. In this
basis, the explicit commutation rules for $su(3)$ are ($i,j,k=1,2,3$):
\be\begin{array}{l}\label{commutatorsA}
[H_i,H_j]=0, \\[0.3cm]
[H_i,F_{jk}]=\hbar\,(\delta_{ij} - \delta_{ik})F_{jk}, \\[0.3cm]
[F_{ij},F_{kl}]=\hbar\,(\delta_{jk} F_{il}- \delta_{il} F_{kj}) +
\hbar\,\delta_{jk} \delta_{il}(H_i - H_j) .
\end{array}\ee

The canonical Lie bialgebra structure  is determined by the
cocommutator: \be\begin{array}{l}\label{cocommutadorA}
\delta (H_i)=0,\\[0.3cm]
\delta (F_{ij})= \frac{1}{2}\; (H_i-H_j) \wedge F_{ij} +\,
\sum_{k=i+1}^{j-1}{F_{ik} \wedge F_{kj}}
\qquad\qquad (i<j),\\[0.3cm]
\delta (F_{ij})=  \frac{1}{2}\; (H_j-H_i) \wedge F_{ij} - \,
\sum_{k=j+1}^{i-1}{F_{ik} \wedge F_{kj}}   \qquad\qquad (i>j).
\end{array}\ee


\subsection{The quantum algebra $U_{z,\hbar}(su(3))$}

By making use of the quantization approach described in \cite{ballesterosJPA08} onto the canonical Lie bialgebra structure  (\ref{commutatorsA})-(\ref{cocommutadorA}), a long but straightforward computation  gives rise to the following coproduct:  
\bea
 \Delta_{z,\hbar} (H_i) &=& H_i \otimes 1 + 1 \otimes H_i
\nonumber\\[0.2cm]
\Delta_{z,\hbar} (F_{12}) &=& e^{\z\,(H_1-H_2)/2} \otimes F_{12} + F_{12} \otimes
e^{-\z\,(H_1-H_2)/2}
\nonumber\\[0.2cm]
\Delta_{z,\hbar} (F_{23}) &=& e^{\z\,(H_2-H_3)/2} \otimes F_{23} + F_{23}
\otimes e^{-\z\,(H_2-H_3)/2}
\nonumber\\[0.2cm]
 \Delta_{z,\hbar} (F_{21}) &=& e^{\z\,(H_1-H_2)/2} \otimes F_{21} + F_{21} \otimes
e^{-\z\,(H_1-H_2)/2}
\nonumber\\[0.2cm]
 \Delta_{z,\hbar} (F_{32}) &=& e^{\z\,(H_2-H_3)/2} \otimes F_{32} + F_{32} \otimes
e^{-\z\,(H_2-H_3)/2}
\nonumber\\[0.2cm]
\Delta_{z,\hbar}(F_{13})&=& \; \,
 e^{\z\,(H_1-H_3)/2} \otimes F_{13} + F_{13} \otimes
e^{-\z\,(H_1-H_3)/2} 
\nonumber\\[0.2cm]
&& \qquad
+ \frac{2}{\hbar} \,\sinh (z\,\hbar/2)\, ( e^{\z\,(H_2-H_3)/2}\,F_{12}
 \otimes e^{-\z\,(H_1-H_2)/2}\,F_{23} 
\nonumber\\[0.2cm]
&& \qquad\qquad -
\, e^{\z\,(H_1-H_2)/2}\,F_{23} \otimes
e^{-\z\,(H_2-H_3)/2}\,F_{12})
\nonumber\\[0.2cm]
 \Delta_{z,\hbar}(F_{31}) &=& 
 e^{\z\,(H_1-H_3)/2} \otimes F_{31} + F_{31} \otimes
e^{-\z\,(H_1-H_3)/2} 
\nonumber\\[0.2cm]
&&  \qquad + \frac{2}{\hbar} \,\sinh (z\,\hbar/2)\, (\,
e^{\z\,(H_2-H_3)/2}\,F_{21}
 \otimes e^{-\z\,(H_1-H_2)/2}\,F_{32} 
\nonumber\\[0.2cm]
&& \qquad\qquad - 
\, e^{\z\,(H_1-H_2)/2}\,F_{32} \otimes
e^{-\z\,(H_2-H_3)/2}\,F_{21}). 
\nonumber
\eea
Concerning the deformed commutation rules, the ones in which the elements of the Cartan subalgebra are involved will remain undeformed with respect to (\ref{commutatorsA}). The remaining ones are found to be:
\[\begin{array}{lll}
&  [F_{12},F_{23}] = \hbar\,F_{13}\qquad\qquad
 &[F_{32},F_{21}] = \hbar\,F_{31}
\\[0.3cm]
& \displaystyle\conm{F_{12}}{F_{13}}=\frac{4}{\hbar}\, (\sinh \frac{z\hbar}{2})^2\, F_{12}\,F_{23}\,F_{12} \qquad\qquad
&\displaystyle \conm{F_{13}}{F_{23}}=\frac{4}{\hbar}\, (\sinh \frac{z\hbar}{2})^2 \,  F_{23}\,F_{12}\,F_{23}
\\[0.3cm]
& \displaystyle\conm{F_{31}}{F_{21}}= \frac{4}{\hbar}\, (\sinh \frac{z\hbar}{2})^2\, F_{21}\,F_{32}\,F_{21} \qquad\qquad
& \displaystyle\conm{F_{32}}{F_{31}}= \frac{4}{\hbar}\, (\sinh \frac{z\hbar}{2})^2\,  F_{32}\,F_{21}\,F_{32}
\\[0.3cm]
&  \conm{F_{23}}{F_{21}}=0 \qquad \qquad
&\conm{F_{12}}{F_{32}}=0
\\[0.3cm]
& \displaystyle [F_{12},F_{21}]= \hbar\,\frac{\sinh(\z\,(H_1-H_2))}{z}\qquad\qquad
&\displaystyle [F_{23},F_{32}]=\hbar\,\frac{\sinh(\z\,(H_2-H_3))}{z}
\end{array}\]
\[\begin{array}{l}
\displaystyle  [F_{13},F_{21}]=-\frac{2}{z}
\,\sinh\frac{z\hbar}{2}\,\cosh(\z\,(H_1-H_2+\frac{\hbar}{2}))\,\,F_{23}
\\[0.3cm]
\displaystyle [F_{13},F_{32}]=\frac{2}{z} \,\sinh\frac{z\hbar}{2}\,\cosh(\z\, (H_2-H_3+\frac{\hbar}{2}))\,F_{12}
\\[0.3cm]
\displaystyle [F_{12},F_{31}]= -\frac{2}{z} \,\sinh\frac{z\hbar}{2}
\,\cosh(\z\,(H_1-H_2
-\frac{\hbar}{2}))\,F_{32}
\\[0.3cm]
\displaystyle [F_{23},F_{31}]=\frac{2}{z} \,\sinh\frac{z\hbar}{2}\,\cosh(\z\, (H_2-H_3-\frac{\hbar}{2}))\,F_{21}
\\[0.3cm]
\displaystyle  [F_{13},F_{31}]=\hbar\displaystyle{ \frac{\sinh(\z\,(H_1-H_3))}{z} }
 +\frac{2}{z\hbar} (\sinh\frac{z\hbar}{2})^2 \,\sinh(\z\,(H_1-H_2))\,
 \anticonm{F_{23}}{F_{32}} 
\nonumber\\[0.2cm]
  \displaystyle\qquad\qquad\qquad\qquad
   +\;\frac{2}{z\hbar}
(\sinh\frac{z\hbar}{2})^2\,\sinh(\z\,(H_2-H_3))\,
\anticonm{F_{12}}{F_{21}}
\end{array}
\]


\subsection{The Poisson-Hopf limit}

If we compute the limits (\ref{poislimit}) and (\ref{colimit}) of the above expressions
we get the following Poisson-Hopf algebra $su_q^{\cal P}(3)$ with coproduct
\bea
  \Delta_z^{\cal P} (H_i) &=& H_i \otimes 1 + 1 \otimes H_i, \quad i=1,2,3
\nonumber\\[0.2cm]
 \Delta_z^{\cal P} (F_{12}) &=& e^{\z\,(H_1-H_2)/2} \otimes F_{12} + F_{12} \otimes
e^{-\z\,(H_1-H_2)/2}
\nonumber\\[0.2cm]
 \Delta_z^{\cal P} (F_{23}) &=& e^{\z\,(H_2-H_3)/2} \otimes F_{23} + F_{23}
\otimes e^{-\z\,(H_2-H_3)/2}
\nonumber\\[0.2cm]
   \Delta_z^{\cal P} (F_{21}) &=& e^{\z\,(H_1-H_2)/2} \otimes F_{21} + F_{21} \otimes
e^{-\z\,(H_1-H_2)/2}
\nonumber\\[0.2cm]
 \Delta_z^{\cal P} (F_{32}) &=& e^{\z\,(H_2-H_3)/2} \otimes F_{32} + F_{32} \otimes
e^{-\z\,(H_2-H_3)/2}
\nonumber\\[0.2cm]
 \Delta_z^{\cal P}(F_{13})&=& \; \,
 e^{\z\,(H_1-H_3)/2} \otimes F_{13} + F_{13} \otimes
e^{-\z\,(H_1-H_3)/2} 
\nonumber\\[0.2cm]
&& \qquad
+ z\, ( e^{\z\,(H_2-H_3)/2}\,F_{12}
 \otimes e^{-\z\,(H_1-H_2)/2}\,F_{23} 
\nonumber\\[0.2cm]
&& \qquad\qquad -
\, e^{\z\,(H_1-H_2)/2}\,F_{23} \otimes
e^{-\z\,(H_2-H_3)/2}\,F_{12})
\nonumber\\[0.2cm]
 \Delta_z^{\cal P}(F_{31}) &=& 
 e^{\z\,(H_1-H_3)/2} \otimes F_{31} + F_{31} \otimes
e^{-\z\,(H_1-H_3)/2} 
\nonumber\\[0.2cm]
&&  \qquad + z\, (\,
e^{\z\,(H_2-H_3)/2}\,F_{21}
 \otimes e^{-\z\,(H_1-H_2)/2}\,F_{32} 
\nonumber\\[0.2cm]
&& \qquad\qquad - 
\, e^{\z\,(H_1-H_2)/2}\,F_{32} \otimes
e^{-\z\,(H_2-H_3)/2}\,F_{21}) 
\nonumber
\eea
The deformed Poisson brackets are:
\[\begin{array}{lll}\label{poissoncommutatorsA}
&\{ H_i,H_j\}=0, \qquad &
\{H_i,F_{jk}\}=\hbar\,(\delta_{ij} - \delta_{ik})F_{jk}, \\[0.3cm]
& \pois{F_{12}}{F_{23}}= F_{13}
\qquad & \pois{F_{32}}{F_{21}}=F_{31}
\\[0.3cm]
& \pois{F_{12}}{F_{13}}=\pois{F_{12}}{\pois{F_{12}}{F_{23}}}=z^2\, F_{12}^2\,F_{23}
\qquad & \pois{F_{13}}{F_{23}}=\pois{\pois{F_{12}}{F_{23}}}{F_{23}}=z^2 \,  F_{23}^2\,F_{12}
\\[0.3cm]
& \pois{F_{31}}{F_{21}}=
\pois{\pois{F_{32}}{F_{21}}}{F_{21}}= z^2\, F_{21}^2\,F_{32}
\qquad
& \pois{F_{32}}{F_{31}}=
\pois{F_{32}}{\pois{F_{32}}{F_{21}}}= z^2 \,  F_{32}^2\,F_{21}
\\[0.3cm]
& \pois{F_{23}}{F_{21}}=0 \qquad 
&\pois{F_{12}}{F_{32}}=0.
\\[0.3cm]
& \pois{F_{12}}{F_{21}} = \frac{1}{z}\, \sinh(\z\,(H_1-H_2))
\qquad
& \pois{F_{23}}{F_{32}} =\frac{1}{z}\,\sinh(\z\,(H_2-H_3))
\\[0.3cm]
 & \pois{F_{13}}{F_{21}}=-\,\cosh(\z\,(H_1-H_2))\,\,F_{23} 
\qquad 
& \pois{F_{13}}{F_{32}}=\,\cosh(\z\, (H_2-H_3))\,F_{12}
\\[0.3cm]
 & \pois{F_{12}}{F_{31}}= -
\,\cosh(\z\,(H_1-H_2))\,F_{32}
\qquad
& \pois{F_{23}}{F_{31}}=\,\cosh(\z\,
(H_2-H_3))\,F_{21}
\end{array}
\] 
\[ 
\quad\,\, \pois{F_{13}}{F_{31}}= \frac{\sinh(\z\,(H_1-H_3))}{z}
 +z \,\sinh(\z\,(H_1-H_2))\,
 \,{F_{23}}\,{F_{32}}+ z\,\sinh(\z\,(H_2-H_3))\,
\,{F_{12}}\,{F_{21}} 
\]
The fact that $\Delta_z^{\cal P} $ is a Poisson algebra homomorphism between $su_q^{\cal P}(3)$ and $su_q^{\cal P}(3)\otimes su_q^{\cal P}(3)$ can be proven by direct computation.

\subsection{Casimir functions}

In the case of the quantum algebra $U_{z,\hbar}(su(3))$ the problem of finding its $q$-deformed Casimir operators is a quite difficult one (see for instance the construction for the $U_{q}(sl(3))$ ones given in \cite{Plaza}). Nevertheless, the corresponding central functions in the case of $su_q^{\cal P}(3)$ can 
can be obtained by direct computation. The explicit form for the deformed `second order' Casimir function is
\bea
&&\!\!\!\!\!\!\!\!\! C^{\cal P}_z=\frac{1}{z^2}\left(
\sinh^2\frac{z(H_1 + H_2 - 2\,H_3)}{3} +
 \sinh^2 \frac{z(H_1 + H_3 - 2\,H_2)}{3} +
 \sinh^2\frac{z(H_2 + H_3 - 2\,H_1)}{3}\right)
\nonumber\\[0.2cm]
&&
 +\, 2\,F_{12}\,F_{21}\,\cosh\frac{z(H_1 + H_2 - 2\,H_3)}{3}+
2\,F_{23}\,F_{32}\,\cosh\frac{z(H_2 + H_3- 2\,H_1)}{3}
\nonumber\\[0.2cm]
&&
+ \, 2\,F_{31}\,F_{13}\,\cosh\frac{z(H_3 + H_1- 2\,H_2)}{3} 
+
2\,z\, (F_{12}\,F_{23}\,F_{31} + F_{21}\,F_{32}\,F_{13})\,\sinh \frac{z(H_3 + H_1 - 2\,H_2)}{3}
\nonumber\\[0.2cm]
&&
 +\, 
2\,z^2\,F_{12}\,F_{21}\,F_{32}\,F_{23}\,\cosh\frac{z(H_3 + H_1 - 2\,H_2)}{3}
\nonumber
\eea
The `third order' one reads:
\bea
&&\!\!\!\!\!\!\!\!\! D^{\cal P}_z=-\frac{1}{z^3}
\sinh \frac{z(H_1 + H_2 - 2\,H_3)}{3}\,\sinh \frac{z(H_1 + H_3 - 2\,H_2)}{3}
\,\sinh \frac{z(H_2 + H_3 - 2\,H_1)}{3}
\nonumber\\[0.2cm]
&&
+\, \frac{1}{z}\,F_{12}\,F_{21}\,
\sinh \frac{z(H_1 + H_2 - 2\,H_3)}{3}
+
\frac{1}{z}\,F_{23}\,F_{32}\,
\sinh \frac{z(H_2 + H_3 - 2\,H_1)}{3}
\nonumber\\[0.2cm]
&&
 +\,
\frac{1}{z}\,F_{31}\,F_{13}\,
\sinh\frac{z(H_3 + H_1 - 2\,H_2)}{3}+
\left(F_{12}\,F_{23}\,F_{31} + F_{21}\,F_{32}\,F_{13}
\right)\,\cosh \frac{z(H_3 + H_1 - 2\,H_2)}{3}
\nonumber\\[0.2cm]
&&
+\,
z\,F_{12}\,F_{21}\,F_{32}\,F_{23}\,\sinh \frac{z(H_3 + H_1 - 2\,H_2)}{3}
\nonumber
\eea
These two expressions can be useful in two different directions. Firstly, as the building blocks for the construction of integrable deformations of $su(3)$ classical spin chains through the formalism given in \cite{BR}. Secondly, as a first-order guide to construct the $q$-deformed Casimir operators for the quantum algebra $U_{z,\hbar}(su(3))$ having in mind the analysis of the corresponding $q$-deformed mass formulae \cite{Hern}. 

\normalsize


\section{Conclusions}

The algebraic structures presented in this paper illustrate the full  `hierarchy of complexity' that Hopf algebras provide: quantum algebras $U_{z,\hbar}(g)$ would be the richest and most complex structures (both non-commutative since $\hbar\neq 0$ and non-cocommutative since $z\neq 0$), the $q$-Poisson algebras would be  somehow intermediate (being commutative and non-cocommutative) and the Lie algebras would be the `simplest' ones (non-commutative and cocommutative). 

 In this context, the procedure to obtain a $q$-Poisson analogue of a given quantum algebra is canonically defined. Firstly, one goes from $U_{z}(g)$ to $U_{z,\hbar}(g)$ 
by constructing analytically the deformation from the very beginning and by taking into account explicitly both parameters  $(\hbar,z)$. At first sight, this two-parameter deformation seems to be irrelevant since
 for any finite value of $\hbar$,  by using the Lie algebra automorphism $X\rightarrow \hbar\,X$ and by transforming the deformation parameter as $z\rightarrow z/\hbar$, the quantum algebra $(U_{z,\hbar}(g),\Delta_{z,\hbar})$ can be converted into $(U_{z}(g),\Delta_{z})$.  However, in order to obtain the Poisson analogue, one has to perform the $\hbar\to 0$ limit that {\em does not commute} with this automorphism and gives rise to the right result. Note that this procedure is quite similar to the well known contraction theory of quantum algebras \cite{celeghini92}.

With respect to the Poisson limit (\ref{poislimit}) and (\ref{colimit}) we would like to stress that it is essential to write the quantum algebra in terms of commutation rules, and {\em not} by making use of $q$-commutators, since the Poisson limit  of the latter is not well-defined. This fact did  not allow in the past the construction of $q$-Poisson analogues for algebras of rank greater than one.
 However, the analytical bases approach presented 
 in~\cite{ballesterosJCS08,ballesterosJPA08}  
 provides a quantization framework based on pure commutators, thus leading to a well-defined Poisson limit for arbitrary quantum algebras. Hence, analytical bases looks to have a privileged connection with the semiclassical limit.
 
On the other hand, quantum deformations of simple Lie algebras  are usually described by means of their simple root generators $X_i$ together with their $q$-Serre  relations. Indeed, this approach simplifies the mathematical scheme as it does not make use of the non simple root generators, even if the latter are necessary in the physical applications due to their role as independent symmetries (see, for instance \cite{Luki}--\cite{Afin}).
However, the $q$-Serre relations can be rewritten in terms of commutators (see \cite{ballesteros06,ballesteros07}) and, after introducing explicitly the $\hbar$ parameter they read
\be\begin{array}{lll}
[X_i,[X_i,X_j]]=4 \sinh^2 \frac {z\hbar}{2}\; X_iX_jX_i &\qquad \mbox{if}\quad a_{ij}=-1 ,\\[0.3cm]
[X_{i},X_{j}]=0 &\qquad \mbox{if}\quad a_{ij}=0 .
\end{array}\ee
Therefore, the corresponding $q$-Poisson-Serre relations in the limit $\hbar \to 0$ can be consistently defined as
\be\begin{array}{lll}
\pois{X_{i}}{\pois{X_{i}}{X_{j}}}=z^2\, X_{i}^2\,X_{j} &\qquad \mbox{if}\quad a_{ij}=-1 ,\\[0.3cm]
\pois{X_{i}}{X_{j}}=0 &\qquad \mbox{if}\quad a_{ij}=0 .
\end{array}\ee


\section*{Acknowledgments}

This work was partially supported  by the Ministerio de
Educaci\'on y Ciencia  of Spain (Projects FIS2005-03989 and MTM2007-67389, both with EU-FEDER support),  by the
Junta de Castilla y Le\'on   (Project GR224), and by
INFN-MICINN (Italy-Spain).





\begin{thebibliography}{99}

\small


\bibitem{KR} Kulish P P and Reshetikhin N 1981 {\it Proc. Steklov Math.
Inst.} {\bf 101}  101

\bibitem{Takh} Takhtajan L A  1990
{\it Introduction to quantum group and integrable massive models of Quantum Field Theory}, Mo-Lin Ge and Bao-Heng Zao (eds.),
(World Scientific, Singapore) p. 69 

\bibitem {FRT} 
  Reshetikhin N Y,  Takhtadzhyan L A and    Faddeev L D 1990  
{\it Leningrad Math. J.} {\bf 1}  193   

\bibitem{CP}  Chari V and  Pressley A 1994 {\em A Guide to Quantum Groups} (Cambridge Univ. Press, Cambridge)

\bibitem{majid}  
Majid S  1995
{\em Foundations of Quantum Group Theory} (Cambridge Univ. Press,
Cambridge)

\bibitem{PS}
Pasquier V and Saleur H 1990 {\it Nucl. Phys.} {\bf B330} 523

\bibitem{KulSkly} Kulish P P and Sklyanin E K 1991 {\it J.
Phys. A: Math. Gen.} {\bf 24} L435


\bibitem{Sierra} 
G\'omez C, Ruiz-Altaba M and Sierra G 1995 {\it
Quantum Groups in Two-Dimensional Physics}, Cambridge University Press 

\bibitem{phonons} Bonechi F,  Celeghini E,  Giachetti R,  Sorace E and  Tarlini M 1992 {\it Phys. Rev. Lett.} {\bf 68}  3718

\bibitem{magnons} Bonechi F,  Celeghini E,  Giachetti R,  Sorace E and  Tarlini M 1992 {\it Phys. Rev.} {\bf B46}  5727

\bibitem{Luki}
 Lukierski J,  Nowicky A and  Ruegg H  1992      
  {\it Phys. Lett.}   {\bf B293},  344

\bibitem{Giller}
Giller S, Koshinski P, Majewski M, Maslanka P and Kunz J 1992      
  {\it Phys. Lett.}   {\bf B286},  57

\bibitem{Ruegg}
Ruegg H  1993
{\it Integrable Systems, Quantum groups and Quantum Field Theories} NATO ASI Series, {\rm L.A. Ibort and M.A. Rodr\'{\i}guez, eds}, p. 45.
  
\bibitem{Dobrev} Dobrev V K 1993 {\it J.
Phys. A: Math. Gen.} {\bf 26} 1317


\bibitem{Afin}
 Ballesteros  A,  Herranz F J , del Olmo M A and  Santander    M     1994   {\it J. Math. Phys.}  {\bf 35},  4928
 
 
 

\bibitem{nuclei} Celeghini E,  Giachetti R,  Sorace E and  Tarlini M 1992 {\it Phys. Lett.} {\bf B280} 180

\bibitem{Bonatsos}
 Bonatsos D and  Daskaloyannis C 1999 {\it Prog. Part. Nucl. Phys}. {\bf 43} 537
 
 
 \bibitem{Gau} Gaudin M 1998 {\it J.
de Physique} {\bf 37} 1087

\bibitem{Skly} Sklyanin E K 1992 {\it Commun. Math. Phys.} {\bf 150} 181

\bibitem{Musso} Falqui G and Musso F 2003 {\it J.
Phys. A: Math. Gen.} {\bf 46} 11655

\bibitem{Lerma} Lerma H S and Errea B F 2003 {\it J.
Phys. A: Math. Theor.} {\bf 2007} 4125

 
 
\bibitem{BR} Ballesteros A and Ragnisco O 1998 {\it J.
Phys. A: Math. Gen.} {\bf 31} 3791

\bibitem{DaP} Da Providencia J 1993 {\it J.
Phys. A: Math. Gen.} {\bf 26} 5845

\bibitem{ballesterosJCS08}  Ballesteros A,   Celeghini E and  del Olmo M A 2008  {\it J. Phys.  C. S.} {\bf 128}  012049

\bibitem{ballesterosJPA08} 
Celeghini E, Ballesteros A    and  del Olmo M A 2008  
{\it J. Phys. A: Math. Theor.} {\bf 41}  304038

\bibitem{Dri87}  Drinfeld V G 1987 ``Quantum Groups'' in 
{\it Proceed. of
the Inter. Congress of Mathematicians}, Berkeley, 1986, A.M. Gleason (ed.)    (AMS, Providence)

\bibitem{Jimbo}  Jimbo M 1985
{\it Lett. Math. Phys}.  {\bf 10}  63   

\bibitem{Ohn} Ohn C  1992 
{\it Lett. Math. Phys.} {\bf 25} 85

\bibitem{ballesteros06}
 Ballesteros A,   Celeghini E and  del Olmo M A 2006
 {\it J. Phys. A: Math. Gen.} {\bf 39}  9161

 \bibitem{ballesteros07}
 Ballesteros A,   Celeghini E and  del Olmo M A 2007
 {\it J. Phys. A: Math. Theor.}  {\bf 40}    2013

\bibitem{Plaza} Rodr\'{\i}guez-Plaza M J  1991 
{\it J. Math. Phys.} {\bf 32} 2020

\bibitem{Hern} C\'arcamo Hern\'andez A E  2006 
{\it J. Math. Phys.} {\bf 47} 073509

\bibitem{celeghini92} 
Celeghini E,  Giachetti R,  Sorace E and  Tarlini M 1992
{\it Lecture Notes in Mathematics} 1510, Kulish P P (ed.),
(Springer, Berlin) p. 221.


\end{thebibliography}
\end{document}